%
\font\fifteenrm=cmr10 scaled\magstep2 
\font\fifteeni=cmmi10 scaled\magstep2
\font\fifteensy=cmsy10 scaled\magstep2
\font\fifteenbf=cmbx10 scaled\magstep2
\font\fifteentt=cmtt10 scaled\magstep2
\font\fifteenit=cmti10 scaled\magstep2
\font\fifteensl=cmsl10 scaled\magstep2
\font\fifteenam=msam10 scaled\magstep2
\font\fifteenbm=msbm10 scaled\magstep2
\font\fifteenex=cmex10 scaled\magstep2
\font\fifteensc=cmcsc10 scaled\magstep2 
\font\twelverm=cmr10 at 12pt
\font\twelvei=cmmi10 at 12pt
\font\twelvesy=cmsy10 at 12pt
\font\twelvebf=cmbx10 at 12pt
\font\twelvett=cmtt10 at 12pt
\font\twelveit=cmti10 at 12pt
\font\twelvesl=cmsl10 at 12pt
\font\twelveam=msam10 at 12pt
\font\twelvebm=msbm10 at 12pt
\font\twelveex=cmex10 at 12pt
\font\twelvesc=cmcsc10 at 12pt
\font\elevenrm=cmr10 scaled\magstephalf 
\font\eleveni=cmmi10 scaled\magstephalf
\font\elevensy=cmsy10 scaled\magstephalf
\font\elevenbf=cmbx10 scaled\magstephalf
\font\eleventt=cmtt10 scaled\magstephalf
\font\elevenit=cmti10 scaled\magstephalf
\font\elevensl=cmsl10 scaled\magstephalf
\font\elevenam=msam10 scaled\magstephalf
\font\elevenbm=msbm10 scaled\magstephalf
\font\elevenex=cmex10 scaled\magstephalf
\font\elevensc=cmcsc10 scaled\magstephalf
\font\tenrm=cmr10
\font\teni=cmmi10
\font\tensy=cmsy10
\font\tenbf=cmbx10
\font\tentt=cmtt10
\font\tenit=cmti10
\font\tensl=cmsl10
\font\tenam=msam10
\font\tenbm=msbm10
\font\tenex=cmex10
\font\tensc=cmcsc10
\font\ninerm=cmr9
\font\ninei=cmmi9
\font\ninesy=cmsy9
\font\ninebf=cmbx9
\font\ninett=cmtt9
\font\nineit=cmti9
\font\ninesl=cmsl9
\font\nineam=msam9
\font\ninebm=msbm9
\font\nineex=cmex9
\font\ninesc=cmcsc9
\font\eightrm=cmr8
\font\eighti=cmmi8
\font\eightsy=cmsy8
\font\eightbf=cmbx8
\font\eighttt=cmtt8
\font\eightit=cmti8
\font\eightsl=cmsl8
\font\eightam=msam8
\font\eightbm=msbm8
\font\eightex=cmex8
\font\eightsc=cmcsc8
\font\sevenrm=cmr7
\font\seveni=cmmi7
\font\sevensy=cmsy7
\font\sevenbf=cmbx7

\font\sevenam=msam7
\font\sevenbm=msbm7

\font\sixrm=cmr6
\font\sixi=cmmi6
\font\sixsy=cmsy6

\font\sixam=msam6
\font\sixbm=msbm6

\font\fiverm=cmr5
\font\fivei=cmmi5
\font\fivesy=cmsy5

\font\fiveam=msam5
\font\fivebm=msbm5

\font\fourrm=cmr5 at 4pt
\font\fouri=cmmi5 at 4pt
\font\foursy=cmsy5 at 4pt

\font\fouram=msam5 at 4pt
\font\fourbm=msbm5 at 4pt

\skewchar\twelvei='177 \skewchar\eleveni='177\skewchar\teni='177
\skewchar\ninei='177 \skewchar\eighti='177\skewchar\seveni='177 
\skewchar\sixi='177 \skewchar\fivei='177 \skewchar\fouri='177
\skewchar\twelvesy='60 \skewchar\elevensy='60 \skewchar\tensy='60
\skewchar\ninesy='60 \skewchar\eightsy='60 \skewchar\sevensy='60 
\skewchar\sixsy='60 \skewchar\fivesy='60 \skewchar\foursy='60
\newfam\itfam
\newfam\slfam
\newfam\bffam
\newfam\ttfam
\newfam\scfam
\newfam\amfam
\newfam\bmfam
\def\eightbig#1{{\hbox{$\left#1\vbox to 6.5pt{}\voidright $}}}
\def\eightBig#1{{\hbox{$\left#1\vbox to 7.5pt{}\voidright $}}}
\def\eightbigg#1{{\hbox{$\left#1\vbox to 10pt{}\voidright $}}}
\def\eightBigg#1{{\hbox{$\left#1\vbox to 13pt{}\voidright $}}}
\def\ninebig#1{{\hbox{$\left#1\vbox to 7.5pt{}\voidright $}}}
\def\nineBig#1{{\hbox{$\left#1\vbox to 8.5pt{}\voidright $}}}
\def\ninebigg#1{{\hbox{$\left#1\vbox to 11.5pt{}\voidright $}}}
\def\nineBigg#1{{\hbox{$\left#1\vbox to 14.5pt{}\voidright $}}}
\def\tenbig#1{{\hbox{$\left#1\vbox to 8.5pt{}\voidright $}}}
\def\tenBig#1{{\hbox{$\left#1\vbox to 9.5pt{}\voidright $}}}
\def\tenbigg#1{{\hbox{$\left#1\vbox to 12.5pt{}\voidright $}}}
\def\tenBigg#1{{\hbox{$\left#1\vbox to 16pt{}\voidright $}}}
\def\elevenbig#1{{\hbox{$\left#1\vbox to 9pt{}\voidright $}}}
\def\elevenBig#1{{\hbox{$\left#1\vbox to 10.5pt{}\voidright $}}}
\def\elevenbigg#1{{\hbox{$\left#1\vbox to 14pt{}\voidright $}}}
\def\elevenBigg#1{{\hbox{$\left#1\vbox to 17.5pt{}\voidright $}}}
\def\twelvebig#1{{\hbox{$\left#1\vbox to 10pt{}\voidright $}}}
\def\twelveBig#1{{\hbox{$\left#1\vbox to 11pt{}\voidright $}}}
\def\twelvebigg#1{{\hbox{$\left#1\vbox to 15pt{}\voidright $}}}
\def\twelveBigg#1{{\hbox{$\left#1\vbox to 19pt{}\voidright $}}}
\def\fifteenbig#1{{\hbox{$\left#1\vbox to 12pt{}\voidright $}}}
\def\fifteenBig#1{{\hbox{$\left#1\vbox to 13.5pt{}\voidright $}}}
\def\fifteenbigg#1{{\hbox{$\left#1\vbox to 18pt{}\voidright $}}}
\def\fifteenBigg#1{{\hbox{$\left#1\vbox to 23pt{}\voidright $}}}
\def\voidright{\right.\nulldelimiterspace=0pt \mathsurround=0pt }
\def\fifteenpoint{
  \textfont0=\fifteenrm \scriptfont0=\twelverm \scriptscriptfont0=\tenrm
  \def\rm{\fam0 \fifteenrm}%
  \textfont1=\fifteeni \scriptfont1=\twelvei \scriptscriptfont1=\teni
  \textfont2=\fifteensy \scriptfont2=\twelvesy \scriptscriptfont2=\tensy
  \textfont3=\fifteenex \scriptfont3=\fifteenex \scriptscriptfont3=\fifteenex
  \def\it{\fam\itfam\fifteenit}\textfont\itfam=\fifteenit
  \def\sl{\fam\slfam\fifteensl}\textfont\slfam=\fifteensl
  \def\bf{\fam\bffam\fifteenbf}\textfont\bffam=\fifteenbf 
    \scriptfont\bffam=\twelvebf\scriptscriptfont\bffam=\tenbf
  \def\tt{\fam\ttfam\fifteentt}\textfont\ttfam=\fifteentt
  \def\sc{\fam\scfam\fifteensc}\textfont\scfam=\fifteensc
  \def\am{\fam\amfam\fifteenam}\textfont\amfam=\fifteenam
    \scriptfont\amfam=\twelveam\scriptscriptfont\amfam=\tenam
  \def\bm{\fam\bmfam\fifteenbm}\textfont\bmfam=\fifteenbm
    \scriptfont\bmfam=\twelvebm\scriptscriptfont\bmfam=\tenbm
  \baselineskip=21pt \rm
  \let\big=\fifteenbig\let\Big=\fifteenBig\let\bigg=\fifteenbigg
  \let\Bigg=\fifteenBigg}
\def\twelvepoint{
  \textfont0=\twelverm \scriptfont0=\ninerm \scriptscriptfont0=\sevenrm
  \def\rm{\fam0 \twelverm}%
  \textfont1=\twelvei \scriptfont1=\ninei \scriptscriptfont1=\seveni
  \textfont2=\twelvesy \scriptfont2=\ninesy \scriptscriptfont2=\sevensy
  \textfont3=\twelveex \scriptfont3=\twelveex \scriptscriptfont3=\twelveex
  \def\it{\fam\itfam\twelveit}\textfont\itfam=\twelveit
  \def\sl{\fam\slfam\twelvesl}\textfont\slfam=\twelvesl
  \def\bf{\fam\bffam\twelvebf}\textfont\bffam=\twelvebf 
    \scriptfont\bffam=\ninebf\scriptscriptfont\bffam=\sevenbf
  \def\tt{\fam\ttfam\twelvett}\textfont\ttfam=\twelvett
  \def\sc{\fam\scfam\twelvesc}\textfont\scfam=\twelvesc
  \def\am{\fam\amfam\twelveam}\textfont\amfam=\twelveam
    \scriptfont\amfam=\nineam\scriptscriptfont\amfam=\sevenam
  \def\bm{\fam\bmfam\twelvebm}\textfont\bmfam=\twelvebm
    \scriptfont\bmfam=\ninebm\scriptscriptfont\bmfam=\sevenbm
  \baselineskip=17.8pt \rm 
  \def\looselineskip{\baselineskip=18.5pt plus 1.8pt}%
  \def\tightlineskip{\baselineskip=16.5pt}%
  \def\verytightlineskip{\baselineskip=15pt}%
  \let\big=\twelvebig\let\Big=\twelveBig\let\bigg=\twelvebigg
  \let\Bigg=\twelveBigg  }
\def\elevenpoint{
  \textfont0=\elevenrm \scriptfont0=\ninerm \scriptscriptfont0=\sixrm
  \def\rm{\fam0 \elevenrm}%
  \textfont1=\eleveni \scriptfont1=\ninei \scriptscriptfont1=\sixi
  \textfont2=\elevensy \scriptfont2=\ninesy \scriptfont2=\sixsy 
  \textfont3=\elevenex \scriptfont3=\elevenex \scriptfont3=\elevenex
  \def\it{\fam\itfam\elevenit}\textfont\itfam=\elevenit
  \def\sl{\fam\slfam\elevensl}\textfont\slfam=\elevensl
  \def\bf{\fam\bffam\elevenbf}\textfont\bffam=\elevenbf
  \def\tt{\fam\ttfam\eleventt}\textfont\ttfam=\eleventt
  \def\sc{\fam\scfam\elevensc}\textfont\scfam=\elevensc
  \def\am{\fam\amfam\elevenam}\textfont\amfam=\elevenam
    \scriptfont\amfam=\nineam\scriptscriptfont\amfam=\sixam
  \def\bm{\fam\bmfam\elevenbm}\textfont\bmfam=\elevenbm
    \scriptfont\bmfam=\ninebm\scriptscriptfont\bmfam=\sixbm
  \baselineskip=15.1pt \rm
  \def\looselineskip{\baselineskip=16pt plus 1.5pt}%
  \def\tightlineskip{\baselineskip=14pt}%
  \def\verytightlineskip{\baselineskip=13pt}%
  \let\big=\elevenbig\let\Big=\elevenBig\let\bigg=\elevenbigg
  \let\Bigg=\elevenBigg  }
\def\tenpoint{
  \textfont0=\tenrm \scriptfont0=\eightrm \scriptscriptfont0=\fiverm
  \def\rm{\fam0 \tenrm}%
  \textfont1=\teni \scriptfont1=\eighti \scriptscriptfont1=\fivei
  \textfont2=\tensy \scriptfont2=\eightsy \scriptfont2=\fivesy 
  \textfont3=\tenex \scriptfont3=\tenex \scriptfont3=\tenex
  \def\it{\fam\itfam\tenit}\textfont\itfam=\tenit
  \def\sl{\fam\slfam\tensl}\textfont\slfam=\tensl
  \def\bf{\fam\bffam\tenbf}\textfont\bffam=\tenbf
  \def\tt{\fam\ttfam\tentt}\textfont\ttfam=\tentt
  \def\sc{\fam\scfam\tensc}\textfont\scfam=\tensc
  \def\am{\fam\amfam\tenam}\textfont\amfam=\tenam
    \scriptfont\amfam=\eightam \scriptscriptfont\amfam=\fiveam
  \def\bm{\fam\bmfam\tenbm}\textfont\bmfam=\tenbm
    \scriptfont\bmfam=\eightbm \scriptscriptfont\bmfam=\fivebm
  \baselineskip=14pt \rm
  \def\looselineskip{\baselineskip=14.8pt plus1.5pt}
  \def\tightlineskip{\baselineskip=13.6pt}%
  \def\verytightlineskip{\baselineskip=13pt}%
  \let\big=\tenbig\let\Big=\tenBig\let\bigg=\tenbigg\let\Bigg=\tenBigg  }
\def\ninepoint{
  \textfont0=\ninerm \scriptfont0=\sevenrm \scriptscriptfont0=\fourrm
  \def\rm{\fam0 \ninerm}%
  \textfont1=\ninei \scriptfont1=\seveni \scriptscriptfont1=\fouri
  \textfont2=\ninesy \scriptfont2=\sevensy \scriptfont2=\foursy 
  \textfont3=\nineex \scriptfont3=\nineex \scriptfont3=\nineex
  \def\it{\fam\itfam\nineit}\textfont\itfam=\nineit
  \def\sl{\fam\slfam\ninesl}\textfont\slfam=\ninesl
  \def\bf{\fam\bffam\ninebf}\textfont\bffam=\ninebf
  \def\tt{\fam\ttfam\ninett}\textfont\ttfam=\ninett
  \def\sc{\fam\scfam\ninesc}\textfont\scfam=\ninesc
  \def\am{\fam\amfam\nineam}\textfont\amfam=\nineam
    \scriptfont\amfam=\nineam\scriptscriptfont\amfam=\fouram
  \def\bm{\fam\bmfam\ninebm}\textfont\bmfam=\ninebm
    \scriptfont\bmfam=\ninebm\scriptscriptfont\bmfam=\fourbm
  \baselineskip=12.6pt \rm
  \let\big=\ninebig\let\Big=\nineBig\let\bigg=\ninebigg
  \let\Bigg=\nineBigg  }
\def\eightpoint{
  \textfont0=\eightrm \scriptfont0=\fiverm \scriptscriptfont0=\fourrm
  \def\rm{\fam0 \eightrm}%
  \textfont1=\eighti \scriptfont1=\fivei \scriptscriptfont1=\fouri
  \textfont2=\eightsy \scriptfont2=\fivesy \scriptfont2=\foursy 
  \textfont3=\eightex \scriptfont3=\eightex \scriptfont3=\eightex
  \def\it{\fam\itfam\eightit}\textfont\itfam=\eightit
  \def\sl{\fam\slfam\eightsl}\textfont\slfam=\eightsl
  \def\bf{\fam\bffam\eightbf}\textfont\bffam=\eightbf
  \def\tt{\fam\ttfam\eighttt}\textfont\ttfam=\eighttt
  \def\sc{\fam\scfam\eightsc}\textfont\scfam=\eightsc
  \def\am{\fam\amfam\eightam}\textfont\amfam=\eightam
    \scriptfont\amfam=\eightam\scriptscriptfont\amfam=\fouram
  \def\bm{\fam\bmfam\eightbm}\textfont\bmfam=\eightbm
    \scriptfont\bmfam=\eightbm\scriptscriptfont\bmfam=\fourbm
  \baselineskip=11.2pt \rm
  \let\big=\eightbig\let\Big=\eightBig\let\bigg=\eightbigg
  \let\Bigg=\eightBigg  }

\twelvepoint
\nopagenumbers
\hsize=6in\vsize=8.8in

\parskip=1pt plus 1pt

\newif\ifSpecialhead\Specialheadfalse
\newbox\specialheadbox

\def\specialhead #1\par{\Specialheadtrue\setbox\specialheadbox=\hbox{#1}}
\headline={{\ifSpecialhead\box\specialheadbox\global\Specialheadfalse\else
     \ifnum\pageno<0{\hfill\quad{\twelvebf\folio}}%
     \else\ifnum\pageno<2\hfill
     \else\hfill\twelvepoint\sc\firstmark\quad{\twelvebf\folio}\fi\fi\fi}}

\def\title#1\par{\bigskip{\def\cr{\par\center}\center\fifteenbf #1\par}\medskip}
\def\subtitle#1\par{\centerline{\fifteenrm #1}\medskip}
\def\author#1\par{\medskip{\def\cr{\par\center\twelvesc}\fifteensc\center#1\par}}
\def\center#1\par{\hfil #1\hfil\par}
\def\abstract.#1\par{\message{Abstract.}%
                    \medskip{\narrower\narrower\tenpoint\tightlineskip
                        \noindent{\bf Abstract.}#1\par}\medskip\noindent}
\def\bigabstract.#1\par{\message{Abstract.}%
                         \medskip{\narrower\narrower\tightlineskip
                         \noindent{\bf Abstract. }#1\par}\medskip\noindent}
\def\acknowledgement#1\par{\footnote{}{#1}}
\def\sectionskip{\Goodbreak\vskip 25pt plus 15pt minus 5pt}
\def\secnumber{\ifquiet
               \else\ifNoSections
                    \else\sectionsymbol\the\secno\quad\fi\fi}
\def\section#1\par{ \NoSectionsfalse\par\sectionskip\proofdepth=0\claimno=0
 \ifquiet\else\advance\secno by1\fi\toks0={#1}
 \immediate\write16{\ifquiet\else Section \the\secno\space\fi
                    \the\toks0}%
 \mark{\secnumber #1}%
 {\fifteenpoint\bf\noindent\secnumber #1}\nobreak\bigskip\quietoff
 \nobreak\noindent}
\def\quiet{\quiettrue}

\def\quietoff{\ifQUIET\else\quietfalse\fi}
\newif\ifquiet
\newif\ifQUIET
\newif\ifNoSections
\newcount\claimtype
\newcount\secno
\newcount\claimno
\newcount\subclaimno
\newcount\subsubclaimno
\newcount\subsubsubclaimno
\newcount\proofdepth
\def\subclaimnumber{\ifquiet\else\ifcase\subclaimno\or A\or B\or C\or D\or E\or
     F\or G\or H\or I\or J\or K\or L\or M\or N\or O\or P\fi\fi}
\def\subsubclaimnumber{\ifquiet\else\ifcase\subsubclaimno\or i\or ii\or iii\or 
   iv\or v\or vi\or vii\or viii\or ix\or x\or xi\or xii\or xiii\or xiv\fi\fi}
\def\subsubsubclaimnumber{\ifquiet\else\ifcase\subsubsubclaimno\or a\or b\or 
   c\or d\or e\or f\or g\or \or h\or i\or j\or k\or l\or m\or n\or o\fi\fi}
\def\claimtag{\ifquiet\else
  \ifNoSections
    \ifcase\proofdepth\the\claimno%
    \or\the\claimno.\subclaimnumber
    \or\the\claimno.\subclaimnumber.\subsubclaimnumber
    \or\the\claimno.\subclaimnumber.\subsubclaimnumber
                                                .\subsubsubclaimnumber\fi
  \else
    \ifcase\proofdepth\the\secno.\the\claimno
    \or\the\secno.\the\claimno.\subclaimnumber
    \or\the\secno.\the\claimno.\subclaimnumber.\subsubclaimnumber
    \or\the\secno.\the\claimno.\subclaimnumber.\subsubclaimnumber
                                                .\subsubsubclaimnumber\fi\fi\fi}
\secno=0\claimno=0\proofdepth=0\subclaimno=0\subsubclaimno=0\subsubsubclaimno=0
\NoSectionstrue
\newbox\qedbox
\def\claimname{\ifcase\claimtype Theorem\or Lemma\or Claim\or Corollary\or
               Question\or Definition\or Remark\or Conjecture\fi}
\def\preclaimskip{\removelastskip
    \ifcase\claimtype\goodbreak\vskip 8pt plus 4pt minus 2pt
                  \or\goodbreak\vskip 6pt plus 4pt minus 1pt
                  \or\goodbreak\vskip 5pt plus 4pt minus 1pt
                  \or\goodbreak\vskip 8pt plus 4pt minus 2pt
                  \or\vskip 7pt plus 4pt minus 2pt
                  \or\vskip 7pt plus 4pt minus 2pt
                  \or\vskip 7pt plus 4pt minus 2pt
                  \or\goodbreak\vskip 8pt plus 4pt minus 2pt\fi}
\def\postclaimskip{\ifcase\claimtype         \vskip 4pt plus 2pt minus 2pt
                                          \or\vskip 3pt plus 2pt minus 2pt
                                          \or\vskip 2pt plus 2pt minus 1pt
                                          \or\vskip 4pt plus 2pt minus 2pt
                                          \or\vskip 1pt plus 2pt 
                                          \or\vskip 4pt plus 4pt 
                                          \or\vskip 3pt plus 2pt
                                          \or\vskip 4pt plus 2pt minus 2pt\fi}
\def\claimfont{\ifcase\claimtype
                  \sl\or\sl\or\sl\or\sl\or\sl\or\rm\or\rm\or\sl\fi}
\def\advancetag{\ifcase\proofdepth\advance\claimno by1
                               \or\advance\subclaimno by1
                               \or\advance\subsubclaimno by1
                               \or\advance\subsubsubclaimno by1\fi}
\def\sayclaim#1.#2 #3\par{\ifquiet\else\advancetag\fi
    \preclaimskip\setbox1=\hbox{#1}\setbox2=\hbox{#2}%
    \toks0={#1 }
    \immediate\write16{\ifdim\wd1>0pt\the\toks0
                       \else\claimname\space\fi \claimtag.}%
    \vbox{\noindent
    {\bf\ifdim\wd1=0pt \claimname\else #1\fi\ifquiet.\else\ \claimtag{\ifNoSections.\fi}\fi}%
    \enspace{\ifdim\wd2>0pt\sc #2\enspace\fi}%
    {\claimfont #3\par}}\postclaimskip\quietoff}
\def\theorem{\claimtype=0\sayclaim}
\def\lemma{\claimtype=1\sayclaim}

\def\question{\claimtype=4\sayclaim}

\def\point#1. #2\par{\item{\rm #1.}#2\par}
\def\points#1\cr\par{\medskip\vbox{\let\cr=\point\point#1\par}\par}
\def\df{\it}
\def\prooffont{}
\def\proofsize{}
\def\proofindent{}
\def\proofskip{\badbreak\ifcase\claimtype    \vskip 3pt plus 2pt minus 2pt
                                          \or\vskip 2pt plus 2pt minus 2pt
                                          \or\vskip 1pt plus 2pt minus 1pt
                                          \or\vskip 3pt plus 2pt minus 2pt
                                          \or\vskip 1pt plus 2pt 
                                          \or\vskip 2pt plus 4pt 
                                          \or\vskip 1pt plus 2pt
                                          \or\vskip 3pt plus 2pt minus 2pt\fi}

\def\Goodbreak{\vskip0pt plus.5in\penalty-1000\vskip0pt plus-.5in}
\def\goodbreak{\penalty-500}
\def\badbreak{\penalty500}
\def\Badbreak{\penalty1000}
\def\proof{\message{proof}\removelastskip\Badbreak\proofskip\begingroup
  \advance\proofdepth by1
  \setbox\qedbox=\hbox{\halmos\raise2pt\hbox{\fiverm\claimname}}%
  \prooffont\proofsize\proofindent\noindent{\bf Proof: }}
\def\proofof#1:{\message{proof}\removelastskip\Badbreak\proofskip\begingroup
  \advance\proofdepth by1
  \setbox\qedbox=\hbox{\halmos\raise2pt\hbox{\fiverm#1}}%
  \prooffont\proofsize\proofindent\noindent{\bf Proof of #1: }}
\def\cite[#1]{[{\tenrm{#1}}]\message{[#1]}}
\edef\ref#1{\expandafter\global\expandafter\edef#1{\noexpand\claimtag}}
\newwrite\notes
\openout\notes=\jobname.notes
\long\def\unexpandedwrite#1#2{\def\finwrite{\write#1}%
   {\aftergroup\finwrite\aftergroup{\sanitize#2\endsanity}}}
\def\sanitize{\futurelet\next\sanswitch}
\let\stoken=\space
\def\sanswitch{\ifx\next\endsanity
  \else\ifcat\noexpand\next\stoken\aftergroup\space\let\next=\eat
   \else\ifcat\noexpand\next\bgroup\aftergroup{\let\next=\eat
    \else\ifcat\noexpand\next\egroup\aftergroup}\let\next=\eat
     \else\let\next=\copytoken\fi\fi\fi\fi \next}
\def\eat{\afterassignment\sanitize \let\next= }
\long\def\copytoken#1{\ifcat\noexpand#1\relax\aftergroup\noexpand
  \else\ifcat\noexpand#1\noexpand~\aftergroup\noexpand\fi\fi
  \aftergroup#1\sanitize}
\def\endsanity\endsanity{}

\def\note#1#2{\hbox to2in{\strut#1\quad\dotfill\quad#2}}
\def\boxit#1{\setbox4=\hbox{\kern1pt#1\kern1pt}
  \hbox{\vrule\vbox{\hrule\kern1pt\box4\kern1pt\hrule}\vrule}}
\def\halmos{\hbox{\am\char'3}} 
\def\qed#1\par{\message{.                                }\setbox1=\hbox{#1}%
  \ifdim\wd1>0pt\setbox\qedbox=\hbox{\halmos\raise2pt\hbox{\fiverm #1}}\fi
  \kern5pt\lower 2pt\hbox{\box\qedbox}\proofskip\goodbreak\endgroup}

\def\sectionsymbol{\S}
\def\k{\kappa}
\def\g{\gamma}

\def\b{\beta}
\def\d{\delta}

\def\l{\lambda}

\def\I1{\mathop{\hbox{\sc i}_1}}

\def\P{{\mathchoice{\hbox{\bm P}}{\hbox{\bm P}}
         {\hbox{\tenbm P}}{\hbox{\sevenbm P}}}}
\def\Q{{\mathchoice{\hbox{\bm Q}}{\hbox{\bm Q}}
         {\hbox{\tenbm Q}}{\hbox{\sevenbm Q}}}}

\def\card#1{\left|#1\right|}

\def\coll{\mathop{\rm coll}}

\def\unifto{\buildrel\lower 7pt\hbox{$\to$}\over\to}

\def\cp{\mathop{\rm cp}\nolimits}

\def\from{\mathbin{\vbox{\baselineskip=3pt\lineskiplimit=0pt
                         \hbox{.}\hbox{.}\hbox{.}}}}

\def\ZFC{\hbox{\sc zfc}}

\def\plus{^{\scriptscriptstyle +}}

\def\in{\mathrel{\mathchoice{\raise 
1pt\hbox{$\scriptstyle\cal\char'62$}}
         {\raise 1pt\hbox{$\scriptstyle\cal\char'62$}}
         {\raise .5pt\hbox{$\scriptscriptstyle\cal\char'62$}}
         {\hbox{$\scriptscriptstyle\cal\char'62$}}}\penalty700{}}
\def\ni{\mathrel{\mathchoice{\raise 1pt\hbox{$\scriptstyle\cal\char'63$}}
                   {\raise 1pt\hbox{$\scriptstyle\cal\char'63$}}
                   {\raise .5pt\hbox{$\scriptscriptstyle\cal\char'63$}}
                   {\hbox{$\scriptscriptstyle\cal\char'63$}}}\penalty700}
\def\of{\mathrel{\mathchoice{\raise 1pt\hbox{$\scriptstyle\subseteq$}}
                   {\raise 1pt\hbox{$\scriptstyle\subseteq$}}
                   {\raise .5pt\hbox{$\scriptscriptstyle\subseteq$}}
                   {\hbox{$\scriptscriptstyle\subseteq$}}}}
\def\fo{\mathrel{\mathchoice{\raise 1pt\hbox{$\scriptstyle\supseteq$}}
                   {\raise 1pt\hbox{$\scriptstyle\supseteq$}}
                   {\raise .5pt\hbox{$\scriptscriptstyle\supseteq$}}
                   {\hbox{$\scriptscriptstyle\supseteq$}}}}
\def\notin{\mathrel{\mathchoice
  {\raise 1pt\hbox{\rlap{$\scriptstyle\;|$}$\scriptstyle\cal\char'62$}}
  {\raise 1pt\hbox{\rlap{$\scriptstyle\kern2pt 
          |$}$\scriptstyle\cal\char'62$}}
  {\raise .5pt\hbox{\rlap{$\scriptscriptstyle\, |$}$\scriptscriptstyle
      \cal\char'62$}}
  {\hbox{\rlap{$\scriptscriptstyle\, |$}$\scriptscriptstyle
     \cal\char'62$}}}%
  \penalty700}

\def\and{\mathrel{\kern1pt\&\kern1pt}}
\def\iff{\mathrel{\leftrightarrow}}

\def\add{\mathop{\rm Add}\nolimits}

\def\lt{\mathrel{\mathchoice{\scriptstyle<}{\scriptstyle<}
   {\scriptscriptstyle<}{\scriptscriptstyle<}}}
\def\lte{\mathrel{\scriptstyle\leq}}

\def\[#1]{\left[\vphantom{\bigm|}#1\right]}
\def\<#1>{\langle\,#1\,\rangle}

\def\image{\mathbin{\hbox{\tt\char'42}}}
\def\restrict{\mathbin{\hbox{\am\char'26}}}
\def\force{\mathbin{\hbox{\am\char'15}}}

\def\st{\mid}
\def\seq<#1>{{\def\st{\mid\penalty650}\left<\,#1\,\right>}}

\def\set#1{\{\,#1\,\}}

\def\th{{\hbox{\fiverm th}}}

\def\lttheta{{\raise 1pt\hbox{$\scriptstyle<$}\theta}}

\def\I1{\mathop{\hbox{\sc i}_1}}
\def\ltk{{{\scriptstyle<}\k}}
\def\ltl{{{\scriptstyle<}\l}}
\def\ltg{{{\scriptstyle<}\g}}
\def\ltd{{{\scriptstyle<}\d}}

\def\lted{{{\scriptstyle\leq}\d}}

\def\ltel{{{\scriptstyle\leq}\l}}

\def\Ptail{\P_{\fiverm \!tail}}

\def\ltg{{\scriptscriptstyle<}\g}

\def\Gtail{G_{\fiverm tail}}

\font\arrow=line10 scaled \magstep1
\def\makeline#1.{\hbox{\arrow\char#1}}
\def\makearrow#1.#2.{\hbox{\arrow\char#1\llap{\char#2}}}
\def\definelinesandarrows#1.#2.#3.#4.#5.{
   \expandafter\edef\csname#4line\endcsname{\makeline#1.}
   \expandafter\edef\csname#4arrow\endcsname{\makearrow#1.#2.}
   \expandafter\edef\csname#5line\endcsname{\makeline#1.}
   \expandafter\edef\csname#5arrow\endcsname{\makearrow#1.#3.}}
\definelinesandarrows 0.18.9.ne.sw.
\definelinesandarrows 1.21.11.nnne.sssw.
\definelinesandarrows 2.14.13.nnnne.ssssw.
\definelinesandarrows 3.23.15.nnnnne.sssssw.
\definelinesandarrows 4.23.15.nnnnnne.ssssssw.
\definelinesandarrows 10.30.29.nne.ssw.
\definelinesandarrows 16.49.41.neeeeee.swwwwww.
\definelinesandarrows 17.51.43.neeee.swwww.
\definelinesandarrows 19.55.47.nehuh.swhuh.
\definelinesandarrows 24.58.41.neeeeeee.swwwwwww.
\definelinesandarrows 26.62.9.neee.swww.
\definelinesandarrows 33.49.25.neeeee.swwwww.
\definelinesandarrows 35.62.61.nee.sww.
\definelinesandarrows 64.82.73.se.nw.
\definelinesandarrows 65.85.75.ssse.nnnw.
\definelinesandarrows 66.78.77.sssse.nnnnw.
\definelinesandarrows 67.87.79.ssssse.nnnnnw.
\definelinesandarrows 68.87.79.sssssse.nnnnnnw.
\definelinesandarrows 74.94.93.sse.nnw.
\definelinesandarrows 80.113.105.seeeeee.nwwwwww.
\definelinesandarrows 81.115.107.seeee.nwwww.
\definelinesandarrows 99.126.125.see.nww.
\def\sejoin#1#2{\setbox1=\hbox{#1}\setbox2=\hbox{#2}%
  \hbox{\vbox{\hbox{\copy1\kern\wd2}\nointerlineskip
              \hbox{\kern\wd1\box2}}}}
\def\nejoin#1#2{\setbox1=\hbox{#1}\setbox2=\hbox{#2}%
  \hbox{\vbox{\hbox{\kern\wd1\copy2}\nointerlineskip\hbox{\copy1\kern\wd2}}}}
\newdimen\hnudge
\newdimen\vnudge
\newdimen\hnudgedefault
\newdimen\vnudgedefault

\def\SEdefaultnudge{\hnudge=-16pt\vnudge=20pt}
\def\Edefaultnudge{\hnudge=-25pt\vnudge=6pt}
\def\Sdefaultnudge{\hnudge=-8pt\vnudge=20pt}
\def\longEdefaultnudge{\hnudge=-5pt\vnudge=6pt}
\def\nudgeright#1pt{\advance\hnudge by#1pt}
\def\nudgeleft#1pt{\advance\hnudge by-#1pt}
\def\nudgeup#1pt{\advance\vnudge by#1pt}
\def\nudgedown#1pt{\advance\vnudge by-#1pt}
\def\label#1{\smash{\llap{\kern\hnudge
                   \raise\vnudge\rlap{$\scriptstyle#1$}\hfill}}}

\def\SEarrow{\SEdefaultnudge
             \sejoin\seeline{\sejoin\seeline{\sejoin\seeline\seearrow}}}

\def\Sarrow{\Sdefaultnudge\setbox1=\hbox{\SEarrow}
           \hbox{\hskip 10pt\vrule height\ht1\hbox{\arrow\char'77}}}
\def\Earrow{\Edefaultnudge\setbox1=\hbox{\SEarrow}
 \hbox{\raise 2pt\hbox{\vrule height-.4pt depth.8ptwidth\wd1\kern2pt
       \llap{\arrow\char'55}}}}
\def\longEarrow{\longEdefaultnudge\setbox1=\hbox{\SEarrow}
      \rlap{\hskip-1.25\wd1\raise 2pt
            \hbox{\vrule height-.4pt depth.8ptwidth2.5\wd1\kern2pt
            \llap{\arrow\char'55}}}}
\def\trianglediagram#1#2#3#4#5#6{%
    {\def\normalbaselines{\baselineskip0pt\lineskip8pt\lineskiplimit0pt}%
       \matrix{#1& &\cr
               \Sarrow\label{#2}&\SEarrow\label{#3}&\cr
               #4&\Earrow\label{#5}&#6\cr}}}

\looselineskip

\title Universal Indestructibility

\author  Arthur W. Apter${}^*$\cr
            Baruch College\cr
	The City University of New York\cr
     	{\tentt awabb@cunyvm.cuny.edu}\cr

\author Joel David Hamkins${}^*$\cr
	Kobe University and\cr
	The City University of New York\cr
	{\tentt hamkins@postbox.csi.cuny.edu}\cr
	{\tentt http://scholar.library.csi.cuny.edu/users/hamkins}\cr

\abstract. From a suitable large cardinal hypothesis, we provide a model with a supercompact cardinal in which {\df universal indestructibility} holds: every supercompact and partially supercompact cardinal $\k$ is fully indestructible by $\ltk$-directed closed forcing. Such a state of affairs is impossible with two supercompact cardinals or even with a cardinal which is supercompact beyond a measurable cardinal.

\footnote{}{${}^*$Our research has been supported in part by PSC-CUNY grants and by a Collaborative Incentive Grant from the CUNY Research Foundation.}%
Laver's intriguing preparation \cite[Lav78] makes any supercompact cardinal $\k$ indestructible by $\ltk$-directed closed forcing. In his model, however, cardinals which are only partially supercompact are not generally indestructible; indeed, almost all of them are highly destructible. But this needn't be so. We aim in this paper to provide a model of a supercompact cardinal with {\df universal indestructibility}, one in which every supercompact and partially supercompact cardinal $\g$ is fully indestructible by $\ltg$-directed closed forcing. 

\quiet\theorem Main Theorem. If there is a high-jump cardinal, then there is a transitive model with a supercompact cardinal in which universal indestructibility holds.

The high jump cardinals, defined below, have a consistency strength above a supercompact cardinal and below an almost huge cardinal. Modified versions of the Main Theorem will provide a model of a strongly compact cardinal in which universal indestructibility holds for strong compactness and a model of a strong cardinal in which universal indestructibility holds for strongness. Let us begin by proving that the Laver preparation itself does not achieve universal indestructibility.

\theorem Observation. After the Laver preparation, above the first non-trivial stage of forcing every partially supercompact non-supercompact cardinal is destructible.

\proof In fact we will show that even the measurability of such cardinals is destructible. Please recall that the Laver preparation of $\k$ is defined relative to the {\df Laver function} $\ell\from\k\to V_\k$, defined inductively so that if $\g$ is a measurable cardinal and $\ell\image\g\of V_\g$ then $\ell(\g)$ is some $x$ chosen with respect to some fixed well-ordering of $V_\k$ such that for a minimal $\l$ the set $x$ has least rank such that there is no $\l$-supercompact embedding $j:V\to M$ with critical point $\g$ such that $j(\ell\restrict\g)(\g)=x$, if such an $x$ exists. The {\df Laver preparation} of $\k$ is the reverse Easton $\k$-iteration which at stage $\g$ forces with $\ell(\g)$, provided that this is the $\P_\g$-name of a $\ltg$-directed closed poset. Laver \cite[Lav78] proved that this forcing makes the supercompactness of $\k$ indestructible by $\ltk$-directed closed forcing. In fact, every supercompact cardinal $\g\lte\k$ becomes indestructible by $\ltg$-directed closed forcing. In particular, all the supercompact cardinals of the ground model are preserved to the extenstion. 

In \cite[Ham98b] and \cite[Ham$\infty$a] the second author of this paper defined that a forcing notion admits a {\df gap} at $\d$ when it factors as $\P_1*\dot\P_2$ where $\card{\P_1}<\d$ and $\force\dot\P_2$ is $\lted$-strategically closed. It is easy to see, for example, that the Laver preparation admits a gap between any two stages of forcing. The Gap Forcing Theorem of \cite[Ham$\infty$a] shows that after forcing $V[G]$ which admits a gap at some $\d<\g$, any embedding $j:V[G]\to M[j(G)]$ with critical point $\g$ such that $M[j(G)]^\d\of M[j(G)]$ in $V[G]$ is the lift of an embedding in the ground model. Since all ultrapower embeddings and most strongness extender embeddings have this small degree of closure, the theorem shows that no forcing which admits a gap below $\g$ can increase the degree of supercompactness, the degree of strong compactness, and (except possibly for singular limit ordinals) the degree of strongness of $\g$.

Suppose now that $\g$ is a partially supercompact but not supercompact cardinal in $V[G]$, and that $\g$ is above the first nontrivial stage of forcing in the preparation. We claim that $\g$ is destructible in $V[G]$. Since all the supercompact cardinals in $V$ were preserved to $V[G]$ there must be some $\l$ such that $\g$ is not $\l$-supercompact in $V$. Let $g\of\coll(\g,\l)$ be $V[G]$-generic for the forcing to collapse $\l$ to $\g$. If $\g$ remains measurable after this forcing, then since $\g$ and $\l$ now have the same cardinality, $\g$ is also $\l$-supercompact in $V[G][g]$. Since the combined forcing also admits a gap below $\g$, it follows by the Gap Forcing Theorem of \cite[Ham$\infty$a] that $\g$ is $\l$-supercompact in $V$, a contradiction. Thus, the measurability of $\g$ must have been destroyed by $g$, as claimed.\qed Observation

What the proof really shows is that after gap forcing, if a measurable cardinal is indestructible, it must have been supercompact in the ground model. (Indeed, this is exactly Corollary 5.4 of \cite[Ham$\infty$b]). In particular, if we obtain universal indestructibility by gap forcing, all measurable cardinals in the extension began as supercompact cardinals in the ground model. So, since any kind of forcing which resembles the Laver preparation will admit a lot of gaps, we are pushed to make a stronger large cardinal assumption in the ground model.

We will therefore assume the existence of a high-jump cardinal. Such cardinals were used implicitly in \cite[Ham98b].  A cardinal $\k$ is a {\df high-jump} cardinal when there is an embedding $j:V\to M$ with critical point $\k$ such that for some $\theta$ we have $M^\theta\of M$ and $j(f)(\k)<\theta$ for every function $f:\k\to\k$. The next two lemmas give some rough bounds on the consistency strength of a high-jump cardinal. 

\lemma. If $\k$ is almost huge, then $\k$ is the $\k^\th$ high-jump cardinal. 

\proof Let $j:V\to M$ be an almost huge embedding, so that $\cp(j)=\k$ and $M^{\lt j(\k)}\of M$. It is easy to see that $j(\k)$ must be inaccessible in $V$, and consequently the set $\set{j(f)(\k)\st f:\k\to\k}$ is bounded by some $\theta<j(\k)$. Factor the diagram as: $$\trianglediagram{V}{j_0}{j}{M_0}{k}{M}$$where $j_0:V\to M_0$ is the induced $\theta$-supercompact embedding obtained by using $j\image\theta$ as a seed. That is, $j_0$ is the ultrapower by the measure $\mu=\set{X\st j\image\theta\in j(X)}$ and $k$ is the inverse of the collapse of $\set{j(f)(j\image\theta)\st f\in V}$. Since $j_0(f)(\k)\leq k(j_0(f)(\k))=j(f)(\k)<\theta$, the embedding $j_0$ witnesses that $\k$ is a high-jump cardinal in $V$. Since the measure $\mu$ exists also in $M$, it follows that $\k$ is a high-jump cardinal in $M$, and consequently the normal measure induced by $j$ concentrates on high-jump cardinals.\qed

\lemma. If $\k$ is a high-jump cardinal, then $V_\k$ has a proper class of supercompact cardinals. 

\proof Suppose that $\k$ is a high-jump cardinal with witnessing embedding $j:V\to M$. Since $\theta=\sup\set{j(f)(\k)\st f:\k\to\k}$ is a strong limit cardinal, it follows by the $\theta$-closure of $M$ that $\k$ is $\lttheta$-supercompact in $M$. Since by the high-jump property the failure of the degree of supercompactness of $\k$ cannot jump over $\theta$, it follows from this that $\k$ is actually ${\scriptstyle<}j(\k)$-supercompact in $M$. So by reflection $\k$ must be a limit of cardinals which are $\ltk$-supercompact in $V$. In particular, $V_\k$ has a proper class of supercompact cardinals.\qed

Essentially the same argument shows that in $V_\theta$ there is a proper class of supercompact cardinals, and $\k$ is the $\k^\th$ supercompact cardinal. Before proving the Main Theorem, we will need one more simple lemma:

\lemma. If $\g\leq\l$ and the $\l$-supercompactness of $\g$ is destructible by some $\ltg$-directed closed forcing $\Q$, then the $\l$-supercompactness of $\g$ is destructible by some $\ltg$-directed closed forcing $\Q'$ (of perhaps slightly larger cardinality) which leaves no measurable cardinals in the interval $(\l,\card{\Q'}]$. 

\proof If $\Q$ is inadequate, simply let $\Q'=\Q*\dot{\coll}((2^{\l^\ltk})\plus, \card{\Q})$. The collapsing poset cannot revive the $\l$-supercompactness of $\g$, and certainly ensures that there are no measurable cardinals between $\l$ and $\card{\Q'}$. Note that if $\Q$ is inadequate, then $2^{\l^\ltk}<\card{\Q}$, and so $\card{\Q'}$ is at most $\card{\Q}^{{\leq}2^{\l^\ltk}}$.\qed

We are now ready to prove the main theorem:  

\theorem Universal Indestructibility Theorem. If there is a high-jump cardinal, then there is a transitive model of \ZFC\ with a supercompact cardinal in which universal indestructibility holds; every supercompact or partially supercompact cardinal $\d$ is fully indestructible by $\ltd$-directed closed forcing. 

\proof The proof proceeds in a trial by fire. Specifically, assuming $\k$ is a high-jump cardinal in $V$, we will perform a certain reverse Easton forcing iteration $\P$ of length less than $\k$ in which at each stage $\g$ we aim to destroy as much of the supercompactness of $\g$ as is possible with $\ltg$-directed closed forcing. The idea is that any supercompact or partially supercompact cardinal which survives this ordeal is {\it ipso facto} indestructible. The large cardinal hypothesis will guarantee that in fact something does survive and the iteration does not simply destroy everything.

Let's begin the construction. In the usual reverse Easton manner, we will take direct limits at the inaccessible stages and inverse limits otherwise; what remains is to describe the forcing $\Q_\g$ which occurs at each stage $\g$.  Suppose inductively that the iteration $\P_\g$ is defined up to stage $\g$. In the special case that some condition in $\P_\g$ forces that $\g$ is $\ltd$-supercompact in $V[G_\g]$, where $\d$ is the next measurable cardinal above $\g$ in $V[G_\g]$, and moreover the $\ltd$-supercompactness of $\g$ is indestructible over $V[G_\g]$ by $\ltg$-directed closed forcing  of rank less than $\k$, then we stop the construction, declare success and, forcing below this condition, give $V_\d[G_\g]$ as our final desired model. Otherwise, we continue the iteration. In this case, there is some minimal $\eta<\d$ such that the $\eta$-supercompactness of $\g$ is destroyed by some $\ltg$-directed closed $\Q$ of rank below $\k$. By the previous lemma, there are always such posets $\Q$ which leave no measurable cardinals in the interval $(\eta,\card{\Q}]$, and we may assume that $\Q$ has the least possible rank below $\k$. In particular, if $\g$ is not measurable in $V[G_\g]$, then $\Q$ is trivial. By the Axiom of Choice, using for example a fixed well-ordering of $V_\k$, let $\dot\Q_\g$ be a name for some such forcing. Thus, the stage $\g$ forcing $\dot\Q_\g$ destroys as much of the supercompactness of $\g$ as is possible to destroy with the kind of forcing in which we are interested. Furthermore, since nontrivial forcing occurs only at measurable cardinal stages, the next nontrivial stage of forcing will lie beyond both $\eta$ and $\card{\dot\Q_\g}$, and consequently none of the later stages of forcing will ever revive the $\eta$-supercompactness of $\g$.

The trial-by-fire observation is simply that after the stage $\g$ forcing, any degree of supercompactness of $\g$ which survives must in fact be indestructible, for by the minimality of $\eta$ if we could have destroyed more supercompactness we would have. Let's explain this in detail. Suppose $\g$ is $\l$-supercompact in $V[G_{\g+1}]$. Necessarily, $\l<\eta$. Furthermore, no $\ltg$-directed closed forcing $\Q'$ in $V[G_{\g+1}]$ can destroy the $\l$-supercompactness of $\g$, for then the forcing $\Q_\g*\dot\Q'$ would have destroyed the $\l$-supercompactness of $\g$ over $V[G_\g]$, contradicting the minimality of $\eta$. In particular, the $\l$-supercompactness of $\g$ is preserved by the tail forcing $\P_{\g,\b}$ which leads to any of the later models $V[G_\b]$ for $\b<\k$. Furthermore, since the next nontrivial stage of forcing after $\g$ is beyond $\eta$ and $\card{\dot\Q_\g}$, the $\eta$-supercompactness of $\g$ is never restored by the later stages of forcing. Consequently, in all the later models $V[G_\b]$ for $\b<\k$, the partial supercompactness of $\g$ is fully indestructible by $\ltg$-directed closed forcing of rank less than $\k$. 

We therefore claim that if we ever stop the construction and declare success, then we have in fact succeeded. Suppose we stop the construction and declare success at stage $\g<\k$, jumping into the resulting model $V_\d[G_\g]$. Since we declared success, it must be that in this model $\g$ is an indestructible supercompact cardinal, and there are no measurable cardinals above $\g$. And the trial-by-fire observation shows that any partially supercompact cardinal $\bar\g$ below $\g$ becomes indestructible at stage $\bar\g$ and remains so in all the later models, including $V[G_\g]$. Notice that indestructibility by posets of rank less than $\k$ becomes full indestructibility in $V_\d[G_\g]$ because $\d<\k$. Thus, $V_\d[G_\g]$ is a model of a supercompact cardinal with universal indestructibility, as we desired.

To complete the proof, then, it suffices for us to show that indeed at some stage before $\k$ we stop the construction and declare success. Suppose towards a contradiction that we do not, and that $G\of\P=\P_\k$ is $V$-generic for the $\k$-iteration, in which at every stage $\g<\k$ we saw need to continue the iteration. Let $j:V\to M$ be an embedding which witnesses that $\k$ is high-jump supercompact. Thus, for some $\theta$ we have $M^\theta\of M$ and $j(f)(\k)<\theta$ for every function $f:\k\to\k$ in $V$. We may factor the forcing as $j(\P)=\P*\dot\Q*\dot\Ptail$ where $\dot\Q$ is (a name for) the stage $\k$ forcing in $M$ and $\dot\Ptail$ is (a name for) the subsequent forcing up to $j(\k)$. Certainly $G\of\P=j(\P)_\k$ is $M$-generic for the iteration up to stage $\k$ on the $M$-side. For some minimal $\eta$ below the next measurable cardinal $\d$ above $\k$ in $M[G]$, therefore, the forcing $\Q$ has minimal rank below $j(\k)$ such that the $\eta$-supercompactness of $\k$ is destroyed over $M[G]$ and $\Q$ leaves no measurable cardinals in the interval $(\eta,\card{\Q}]$. Since $\eta$, $\d$ and the rank of $\Q$ are easily defined from $\k$, one can easily find functions $f$ so that they are less than $j(f)(\k)$. Thus, by the high-jump property of $\k$, they are all also less than $\theta$. And since $M$ and $V$ agree up to $\theta$, it follows that $\Q$ destroys the $\eta$-supercompactness of $\k$ over $V[G]$ as well. Let $g\of\Q$ be $V[G]$-generic. By further forcing to add $\Gtail\of\Ptail$ over $V[G][g]$ we can lift the embedding to $j:V[G]\to M[j(G)]$ where $j(G)=G*g*\Gtail$. Next, we find a master condition below $j\image g$, force below it to add $j(g)\of j(\Q)$ generically over $V[G][g][\Gtail]$, and lift the embedding to $j:V[G][g]\to M[j(G)][j(g)]$. Let $\mu$ be the $V[G][g]$-measure on $P_\k\eta$ germinated by the seed $j\image\eta$, so that $X\in\mu\iff j\image\eta\in j(X)$. Since the tail forcing $\Ptail*j(\dot\Q)$ is $\ltd$-closed over $M[G][g]$, it is also $\ltd$-closed over $V[G][g]$, and consequently could not have added the measure $\mu$. Thus, $\mu$ must lie in $V[G][g]$, contradicting the assumption that $\Q$ destroyed the $\eta$-supercompactness of $\k$.\qed

The trial-by-fire idea is easily modified to yield the following theorems. 

\theorem Universal Indestructibility Theorem. If there is a high-jump cardinal, then there is a transitive model with a supercompact cardinal in which every supercompact, partially supercompact, measurable, Ramsey and weakly compact cardinal $\g$ is fully indestructible by $\ltg$-directed closed forcing.

\proof In this theorem, we aim at every stage $\g$ first to destroy the weak compactness of $\g$ if possible, and if this is not possible, to destroy the Ramseyness of $\g$, and if this is not possible, to destroy the measurability of $\g$ and if this is not possible, to destroy as much of the supercompactness of $\g$ as possible. And we do so with posets $\Q$ which leave no additional weakly compact cardinals between the failure of the supercompactnss of $\g$ and $\card{\Q}$. The construction terminates when we find a cardinal $\g$ which in $V[G_\g]$ is indestructibly supercompact up to the next weakly compact cardinal $\d$, the resulting model being $V_\d[G_\g]$. The previous argument shows that the construction terminates before the least high-jump cardinal, and the trial-by-fire observation shows that every supercompact, partially supercompact, measurable, Ramsey and weakly compact cardinal in $V_\d[G_\g]$ is fully indestructible, as desired.\qed

\theorem Universal Indestructibility Theorem. If there is a high-jump cardinal, then there is a transitive model with a strongly compact cardinal in which every strongly compact, partially strongly compact, measurable, Ramsey and weakly compact cardinal $\g$ is fully indestructible by $\ltg$-directed closed forcing. 

\proof For this theorem, at stage $\g$ destroy the weak compactness of $\g$ if possible, and failing that, destroy the Ramseyness of $\g$, and so on, until in the end if we cannot destroy the measurability of $\g$ then we destroy as much of the strong compactness of $\g$ as possible. We declare success when we find a cardinal $\g$ which is indestructibly strongly compact up to the next weakly compact cardinal $\d$. As before, this happens before the least high-jump cardinal, and the trial-by-fire observation shows that the surviving cardinals are fully indestructible in the resulting model $V_\d[G_\g]$.\qed

\theorem Universal Indestructibility Theorem. If there is a high-jump cardinal, then there is a transitive model with a strong cardinal in which every strong, partially strong, measurable, Ramsey and weakly compact cardinal $\d$ is fully indestructible by $\ltd$-directed closed forcing. 

\proof For this theorem, we aim at stage $\g$ first to destroy the weak compactness of $\g$, and failing that to destroy the Ramseyness of $\g$, and so on, until if we cannot destroy the measurability of $\g$ then we destroy as much of the strongness of $\g$ as possible. Before the least high-jump cardinal, the construction will find a cardinal $\g$ which in $V[G_\g]$ is indestructibly strong up to the next weakly compact cardinal $\d$, and so the model $V_\d[G_\g]$ is as desired.\qed

The essence of these theorems lies in the trial by fire. The general procedure, given any list of large cardinal properties of the kind mentioned in these theorems, provided they are well-ordered in strength, is to perform a forcing iteration in which at every stage $\g$ we aim to destroy as much of the large cardinal strength of $\g$, with respect to our list of properties, as is possible. The trial-by-fire result is that in the final model, any cardinal having a property on the list will be indestructible, in virtue of having survived. With a suitable large cardinal in the ground model, one then argues that in fact some large cardinals have survived the trial, and a model of the desired sort of universal indestructibility is produced. 

Perhaps the first use of the trial-by-fire technique is the following theorem, proved by the second author in response to the questions of the first author in \cite[Apt98] concerning the possibility of indestructible measurable limits of supercompact cardinals. 

\theorem.(Hamkins) If $V$ is any model of \ZFC, then there is a forcing extension $V[G]$ such that
\points 1. Every isolated supercompact cardinal of $V$ remains supercompact in $V[G]$ and becomes indestructible there.\cr
            2. No new supercompact cardinals are created; indeed, no cardinal has its degree of supercompactness increased from $V$ to $V[G]$.\cr
            3. Every supercompact limit of supercompact cardinals in $V$ remains strongly compact in $V[G]$.\cr
            4. Every measurable limit of supercompact cardinals in $V[G]$ is fully indestructible.\cr

\proof The proof proceeds by folding into the universal Laver preparation some additional forcing which destroys the measurability of the measurable limits of supercompacts, if possible. The result is that all supercompact cardinals in the extension are indestructible and---the essence of the trial by fire---any surviving measurable limit of supercompact cardinals is indestructible, in virtue of having survived. 

Let $\ell$ be the universal Laver (class) function defined as in Observation 1, assuming a class choice principle, which can be easily forced if necessary (see also Lemma 1 of \cite[Apt98]). Let $\P$ be the reverse Easton support class forcing iteration which first adds a Cohen real (in order to introduce a very low gap), and that at stage $\g$ forces with $\ell(\g)$, provided that this is the $\P_\g$-name of a $\ltg$-directed closed forcing notion. In the case that $\g$ happens to be in $V$ a measurable limit of supercompact cardinals whose measurability is destroyed by some $\ltg$-directed closed forcing $\Q$ in $V^{\P_\g}$, then instead we perform some such forcing of least possible rank, chosen by using the least name in $V$ with respect to a fixed class well-ordering of $V$ (and we do so with a poset $\Q$ which leaves no measurable cardinals between $\g$ and $\card{\Q}$). Suppose now that $G\of\P$ is $V$-generic for this forcing. If $\g$ is a supercompact cardinal in $V$, but not a limit of supercompact cardinals, then the usual Laver argument shows that $\g$ becomes an indestructible supercompact cardinal in $V[G]$, so (1) holds. Since the forcing admits a very low gap, the Gap Forcing Theorem of \cite[Ham$\infty$a] implies that (2) holds. The trial-by-fire observation shows that any measurable limit $\k$ of supercompact cardinals which survives to $V[G]$ must in fact be indestructible there, or else we would have destroyed it at stage $\k$, so (4) holds. What remains is to prove that 
nontrivial instances of this can actually occur, that is, that (3) holds. Suppose that $\k$ is a supercompact limit of supercompact cardinals in $V$. In order to show that $\k$ is strongly compact in $V[G]$, it suffices to show that $\k$ is measurable there, since any measurable limit of supercompact cardinals is strongly compact \cite[Men74]. If $\k$ is not measurable in $V[G]$, then by the closure of the tail forcing, it must be that $\k$ is not measurable in $V[G_{\k+1}]$. Let $j:V\to M$ be a $\theta$-supercompactness embedding, where $\theta$ is much larger than the size of the term for the stage $\k$ forcing $\dot\Q$ in $V$. Factor the forcing $j(\P_\k)$ as $\P_\k*\dot{\widetilde\Q}*\dot\Ptail$, where $\dot{\widetilde\Q}$ is a term for the stage $\k$ forcing in $M$. Since $M$ and $V$ agree up to $\theta$, the term $\dot\Q$ exists in $M$ and furthermore forcing with its interpretation $\Q$ destroys the measurability of $\k$ over $M[G]$. Consequently, the rank of $\dot{\widetilde\Q}$ is no greater than that of $\dot\Q$, and is consequently less than $\theta$. Thus, over $V[G]$, the forcing $\widetilde\Q$ destroys the measurability of $\k$. Supposing that $g\of\widetilde\Q$ is $V[G]$-generic, we conclude that $\k$ is not measurable in $V[G][g]$. Force now to add $\Gtail\of\Ptail$ generically over $V[G][g]$ and lift the embedding to $j:V[G]\to M[j(G)]$ 
where $j(G)=G*g*\Gtail$. Now, find a master condition below $j\image g$ in $j(\widetilde\Q)$ and force below it to add $j(g)\of j(\widetilde\Q)$ generically over $V[G][g][\Gtail]$. In the resulting model, we may lift the embedding to $j:V[G][g]\to M[j(G)][j(g)]$. The induced normal measure cannot have been added by the forcing $\Ptail*j(\dot{\widetilde\Q})$, and so lies in $V[G][g]$. Thus, $\k$ is measurable there, contradicting our earlier observation that it was not. So (3) holds.\qed

While one might hope for universal indestructibility in the presence of several supercompact cardinals, the next theorem shows that this is simply inconsistent.

\theorem. If there are two supercompact cardinals, then universal indestructibility fails for partial supercompactness. Indeed, if a cardinal $\k$ is $\l\plus$-supercompact for some measurable cardinal $\l>\k$, then universal indestructibility fails for partial supercompactness. 

\proof Suppose that every measurable cardinal below $\k$ is indestructible and for some measurable cardinal $\l$ above $\k$ the $\l\plus$-supercompactness of $\k$ is indestructible. By further forcing if necessary, we may assume that $2^\l=\l\plus$ because the forcing to achieve this is $\ltel$-directed closed, and therefore preserves the measurability of $\l$ and the indestructibility of the measurable cardinals below $\l$. Suppose now that the $\l\plus$-supercompactness of $\k$ survives the forcing which adds a new Cohen subset $A\of\k$. It is easy to see that the measurable cardinals below $\k$ remain indestructible in $V[A]$. Nevertheless, the Superdestruction Theorem of \cite[Ham98a] shows that the measurability of $\l$ becomes superdestructible in $V[A]$, and is destroyed by any further $\ltl$-closed forcing which adds a subset to $\l$. Suppose now that $j:V[A]\to M[j(A)]$ witnesses the $\l\plus$-supercompactness of $\k$. Since $V[A]$ and $M[j(A)]$ agree up to $\l\plus=2^\l$, they agree that the measurability of $\l$ is destroyed by the forcing $\add(\l,1)^{V[A]}=\add(\l,1)^{M[j(A)]}$. But by elementarity every measurable cardinal below $j(\k)$ is indestructible in $M[j(A)]$, a contradiction.\qed

One might also hope to combine the various kinds of universal indestructibility, and have for example a model with a supercompact cardinal in which partial strongness is indestructible. But this also is inconsistent. Indeed, one cannot even have a cardinal which is strong beyond a measurable cardinal when partial strongness is universally indestructible. 

\theorem. If a cardinal $\k$ is $(\l+2)$-strong for some measurable cardinal $\l>\k$, then universal indestructibility fails for partial strongness. 

\proof Suppose that $\k$ is $(\l+2)$-strong, where $\l>\k$ is measurable, and partial strongness is universally indestructible. Then after forcing to add a new Cohen subset $A\of\k$, there will be an embedding $j:V[A]\to M[j(A)]$ witnessing the $(\l+2)$-strongness of $\k$ in $V[A]$. The cardinal $\l$ will be measurable in $V[A]$ and consequently also in $M[j(A)]$. By the Superdestruction Theorem of \cite[Ham98a], however, the measurability of $\l$ is destroyed over $V[A]$ by any $\ltl$-closed forcing which adds a subset to $\l$. Since this will also be true over $M[j(A)]$, by reflection it follows that there are destructible measurable cardinals below $\k$, contradicting the hypothesis.\qed

One can make a similar argument in the case of weak compactness:

\theorem. If every weakly compact cardinal $\d$ below $\k$ is indestructible by $\add(\d,1)$, then the $(\l+1)$-strongness of $\k$ is destroyed by $\add(\k,1)$ where $\l$ is the next weakly compact cardinal above $\k$. 

\proof Since the weak compactness of $\l$ can be verified in $V_{\l+1}$, the same argument as the previous theorem applies also in this case. The Superdestruction Theorem of \cite[Ham98a] shows that after small forcing, even the weak compactness of $\l$ is destroyed by $\add(\l,1)$.\qed

Let us now turn to the question of the consistency strength of universal indestructibility. Since we proved that the constructions of the Universal Indestructibility Theorems terminate {\it before} the least high-jump cardinal, we naturally expect that the hypothesis can be reduced. One might naively hope to perform a Laver-like preparation and preserve any given supercompact cardinal while making the partial supercompact cardinals below fully indestructible. But as we mentioned just after the initial Observation, the theorem below from \cite[Ham$\infty$b] shows that this is impossible. 

\theorem.(Hamkins) After forcing with a gap below $\k$, if the measurability of $\k$ is indestructible, then $\k$ was supercompact in the ground model.

\noindent Thus, if one hopes to obtain universal indestructibility by forcing with a gap, one must begin with many supercompact cardinals in the ground model. But what is the optimal hypothesis?

\question. What is the consistency strength of a supercompact cardinal in the presence of universal indestructibilty? What about strongly compact cardinals and strong cardinals? 

\question. What is the consistency strength of an indestructible measurable cardinal? What about an indestructible Ramsey cardinal or an indestructible weakly compact cardinal? 

Let us conclude the paper by showing how a modified Laver preparation can ensure a weaker form of universal indestructibility beginning with just one supercompact cardinal. Specifically, we say that universal {\df partial} indestructibility holds when any cardinal $\g$ which is $\l$-supercompact is indestructible by any $\ltg$-directed closed forcing of size at most $\l$. A cardinal $\g$ is supercompact {\df up to} a measurable cardinal when it is $\ltl$-supercompact for some measurable cardinal $\l$. It is easy to see that this implies that $\g$ is also $\ltel$-supercompact.

\theorem Partial Indestructibility Theorem. Suppose that $\k$ is supercompact and no cardinal is supercompact up to a measurable cardinal. Then there is a forcing extension in which $\k$ is an indestructible supercompact cardinal and universal partial indestructibility holds. 

\proof We simply use the trial-by-fire technique while restricting the size of the forcing notions we consider. Specifically, let $\P$ be the reverse Easton $\k$-iteration in which at stage $\g$ the forcing $\Q_\g$ is chosen, with respect to a well-ordering of the names in $V_\k$, so as to destroy as much of the supercompactness of $\g$ over $V[G_\g]$ as is possible with $\card{\Q_\g}<\l_\g$ where $\l_\g$ is least such that $\g$ is not $\l_\g$-supercompact in $V[G_\g]$. If no such forcing exists, then $\Q_\g$ is trivial. In particular, nontrivial forcing occurs only at measurable cardinal stages. Suppose that $G\of\P$ is $V$-generic for this forcing. 

The trial-by-fire observation in this case is that after the stage $\g$ forcing, any remaining supercompactness of $\g$ is indestructible by $\ltg$-directed closed forcing of size less than $\l_\g$, since if we could have destroyed more, we would have. And since the next nontrivial stage of forcing does not occur until the next measurable cardinal, which by assumption is above $\l_\g$, no amount of supercompactness will be revived by the later stages of forcing. Thus, in all the later models, $\g$ will be partially indestructible. 

To finish the proof, therefore, it suffices to show that $\k$ becomes indestructibly supercompact in $V[G]$. Suppose that some $\ltk$-directed closed forcing $\Q$ destroys the $\l$-supercompactness of $\k$. Let $\theta$ be much larger than $\l$ and $\card{\Q}$ and suppose that $j:V\to M$ is a $\theta$-supercompact embedding in the ground model. Since the $\l$-supercompactness of $\k$ is destructible over $V[G]$ by $\Q$, this will also be true in $M[G]$. Thus, the iteration factors as $j(\P)=\P*\dot{\widetilde\Q}*\dot\Ptail$ where $\dot{\widetilde\Q}$ is a name for the stage $\k$ forcing chosen in $M$ which destroys the $\l$-supercompactness of $\k$ (and possibly more) over $M[G]$. Since $M[G]$ and $V[G]$ agree up to $\theta$, it must be that ${\widetilde\Q}$ also destroys the $\l$-supercompactness of $\k$ over $V[G]$. Force over $V[G]$ to add $g\of{\widetilde\Q}$ and also $\Gtail\of\Ptail$ and lift the embedding to $j:V[G]\to M[j(G)]$ where $j(G)=G*g*\Gtail$. Now, using the directed closure of $j({\widetilde\Q})$, find a master condition below $j\image g$ and force below it to add $j(g)\of j({\widetilde\Q})$ generically over $V[G][g][\Gtail]$. In $V[G][g][\Gtail][j(g)]$, the embedding lifts to $j:V[G][g]\to M[j(G)][j(g)]$. The induced $V[G][g]$-measure $\mu$ defined by $X\in\mu\iff j\image\l\in j(X)$ cannot have been added by the forcing $\Ptail*j(\dot{\widetilde\Q})$, and so it must lie in $V[G][g]$, contradicting that ${\widetilde\Q}$ was supposed to destroy the $\l$-supercompactness of $\k$.\qed

\quiet\section Bibliography

\nopagenumbers
\parindent=0pt
\newbox\Article
\newbox\Journal
\newbox\Author
\newbox\Vol
\newbox\No
\newbox\Year
\newbox\Page
\newbox\Book
\newbox\Publisher
\newbox\Pubaddr
\newbox\Key
\newbox\Editor
\newbox\Comment
\newbox\Note
\def\entry#1#2\par{\item{#1\quad}\hskip-1.1em#2\par}
\def\article#1{\setbox\Article=\hbox{\sl #1, }}
\def\journal#1{\setbox\Journal=\hbox{\rm #1 }}
\def\author#1{\setbox\Author=\hbox{\sc #1, }}
\def\vol#1{\setbox\Vol=\hbox{\bf #1 }}
\def\no#1{\setbox\No=\hbox{no. #1 }}
\def\year#1{\setbox\Year=\hbox{\rm({\oldstyle #1}) }}
\def\page#1{\setbox\Page=\hbox{\rm p. #1 }}
\def\book#1{\setbox\Book=\hbox{\it #1, }}
\def\publisher#1{\setbox\Publisher=\hbox{\rm #1, }}
\def\pubaddr#1{\setbox\Pubaddr=\hbox{\rm #1, }}
\def\key#1{\setbox\Key=\hbox{#1}}
\def\editor#1{\setbox\Editor=\hbox{\rm(#1, Ed.) }}
\def\comment#1{\setbox\Comment=\hbox{\rm #1}}
\def\note#1{\setbox\Note=\hbox{\rm #1 }}
\def\ref#1\par{\smallskip{#1
  \entry{\ifhbox\Key\unhbox\Key\else[\ ]\fi}%
  \unhbox\Author\unhbox\Note
  \ifhbox\Book \unhbox\Book\unhbox\Publisher\unhbox\Pubaddr
               \unhbox\Editor\unhbox\Page\unhbox\Year\unhbox\Comment
  \else \unhbox\Article\unhbox\Journal\unhbox\Vol\unhbox\No\unhbox\Editor
        \unhbox\Page\unhbox\Year\unhbox\Comment\fi\par}}

\tenpoint\tightlineskip

\ref
\author{Arthur W. Apter}
\article{Laver indestructibility and the class of compact cardinals}
\journal{Journal of Symbolic Logic}
\year{1998}
\vol{63}
\no{1}
\page{149-157}
\key{[Apt98]}

\ref
\author{Joel David Hamkins}
\article{Small forcing makes any cardinal superdestructible}
\journal{Journal of Symbolic Logic}
\year{1998}
\vol{63}
\no{1}
\page{51-58}
\key{[Ham98a]}

\ref
\author{Joel David Hamkins}
\article{Destruction or preservation as you like it}
\journal{Annals of Pure and Applied Logic}
\year{1998}
\vol{91}
\page{191-229}
\key{[Ham98b]}

\ref
\author{Joel David Hamkins}
\article{Gap forcing}
\journal{submitted to the Bulletin of Symbolic Logic}
\key{[Ham$\infty$a]}

\ref
\author{Joel David Hamkins}
\article{The lottery preparation}
\journal{submitted to the Annals of Pure and Applied Logic}
\key{[Ham$\infty$b]}

\ref
\author{Richard Laver}
\article{Making the supercompactness of $\kappa$ indestructible under 
 $\kappa$-directed closed forcing}
\journal{Israel Journal Math.}
\vol{29}
\year{1978}
\page{385-388}
\key{[Lav78]}

\ref
\author{Telis K. Menas}
\article{On strong compactness and supercompactness}
\journal{Annals of Math. Logic}
\vol{7}
\year{1974}
\page{327-359}
\key{[Men74]}	

\ref
\author{Azriel Levy, Robert M. Solovay}
\article{Measurable cardinals and the Continuum Hypothesis}
\journal{Israel Journal Math.}
\vol{5}
\year{1967}
\page{234-248}
\key{[LevSol67]}

\bye